\documentclass[leqno,11pt,A4]{article}
 \usepackage{amsmath,amsfonts,amssymb,amsthm}
\usepackage{fancyhdr,mathptmx,enumerate,begmacras}
\newcommand{\mr}{\ensuremath{\mathfrak{M}(\bbr)}}

 \title{Entiers al\'eatoires, ensembles de Sidon, densit\'e 
 dans le groupe de Bohr et ensembles d'analyticit\'e}
 \author{Jean-Pierre Kahane et Yitzhak Katznelson\\[5pt]\small
(Note pr\'esent\'ee aux Comptes rendus de l'Acad\'emie des Sciences) 
 }
 \date{} 
\begin{document}
\maketitle

\begin{verse}\qquad
R\'esum\'e.  {\small  On s\'electionne des entiers $n$ au hasard  
ind\'ependamment les uns des autres, avec la probabilit\'e 
$\varpi_{n}$,  et on \'etudie les propri\'et\'es de 
la suite $\Lambda$ obtenue en distinguant 
deux cas : si $\limsup_{n\to\infty}n\varpi_{n}<\infty$, $\Lambda$
est p.s. un ensemble de Sidon non-dense dans le groupe de Bohr ;
si $\lim_{n\to\infty}n\varpi_{n}=\infty$, $\Lambda$ est p.s. dense dans le
groupe de Bohr et n'est pas un ensemble de Sidon, et, de plus, c'est un ensemble d'analyticit\'e.
}
\end{verse}
\b

\begin{center}
{\bf Random sequences of integers, Sidon sets, density in the Bohr group,\\
and sets of analyticity}
\end{center}

\begin{verse}\qquad
Abstract.  {\small  
We study properties of a sequence $\Lambda$ obtained by a random
selection of integers $n$, where $n\in\Lambda$ with probability
$\varpi_{n}$, 
independently of the other choices.
We distinguish two cases  : if $\limsup_{n\to\infty}n\varpi_{n}<\infty$, 
$\Lambda$ is a.s. a Sidon set, non-dense in the Bohr group ;
if $\lim_{n\to\infty}n\varpi_{n}=\infty$, then $\Lambda$ 
is a.s. a set of analyticity and is 
dense in the Bohr group.
}
\end{verse}

\s
\begin{center}
Abridged English version
\end{center}

\s
The work described in this Note deals with 
harmonic-analytic properties  of random sequences
of positive integers : that of being a Sidon set, a set of analyticity,
dense, or non-dense in the Bohr group $\bbb$,
(the Bohr compactification of the  integers,  the dual 
group 
of $\bbt_{d}$,  the circle endowed with the discrete topology.)

A sequence  $\Lambda\subset\bbz$ is a  Sidon set if 
$c_{0}(\Lambda)=A(\Lambda)$, that is :
every  sequence $\{a_{\lambda}\}_{\lambda\in\Lambda}$ 
which tends to $0$ as $\abs{\lambda}\to\infty$
is the restriction to $\Lambda$ 
of $\hat f$ for some  $f\in\lt$.  
$\Lambda$ is a set of analyticity if only analytic
functions operate in $A(\Lambda)$, that is, every function 
$F$ defined on $\bbr$ such that $F\circ \vp\in A(\Lambda)$ for
every real-valued $\vp\in A(\Lambda)$ is analytic 
in some neighborhood of $0$.

Random sequences provide classes of sequences for which there
is a \emph{statistical} answer to some long standing open problems :

\n
\aba The \bfit{dichotomy} problem :  is every sequence $\Lambda\subset\bbz$ 
is either a Sidon set or  a set of analyticity?

\n
\bab Can a Sidon set be dense in $\bbb$?

The class of random sequences studied in this work is different
from those appearing in earlier work, \cite{ykpm},   \cite{yk1}, \cite{yk2}. 
We consider sets obtained from a sequence  $\bw=\{w_{n}\}$,
$0\le w_{n}<1$ by introducing independent 
Poisson variables $\xi_{n}$ with  parameters $w_{n}$, or
Bernoulli zero-one variables 
$\beta_{n}$ with probability $\pr(\beta_{n}=1)=\varpi_{n}=(1-e^{-w_{n}})$
and set $\Lambda=\Lambda_{\bw}(\omega)=\set{n}{\beta_{n}=1}$ (or 
$\Lambda=\set{n}{\xi_{n}\ge1}$). Notice that the sequences 
defined by $\{\xi_{n}\}$
are the same and have the same statistics as those given by $\beta_{n}$, but the Poisson variables
give different weights to some of the selected points, which 
 simplifies some computations.   

For ``reasonably regular'' sequences of parameters we obtain a statistical
dichotomy:  
\begin{theorem1}
If $nw_{n}=\osta{1}$ then  $\Lambda$ is  a.s. a Sidon set and
is non-dense in $\bbb$.
\end{theorem1}
\begin{theorem2}
If $\lim_{n\to\infty}nw_{n}=\infty$ then  $\Lambda$ is  a.s. a  set
of analyticity   dense in $\bbb$.
\end{theorem2}
Observe that ``reasonable regularity'' is needed. The condition 
$\lim_{n\to\infty}nw_{n}=\infty$ in Theorem 2 cannot be replaced
by $\limsup_{n\to\infty}nw_{n}=\infty$ :
if $w_{n}=1/2$ for $n\in\{3^{k}\}_{k=1}^{\infty}$ 
and $w_{n}=1/n$ for all other $n$, the sequence $\Lambda$
obtained is a union of two Sidon sets, and hence Sidon.

\begin{center}{-----------------}\end{center}
\b
 
 Cette note fait suite \`a des travaux anciens,
  \cite{ykpm}, \cite{yk1}, \cite{yk2}, 
 relatifs aux ensembles de Helson et de Sidon d'une part, aux ensembles d'analyticit\'e et \`a la
 densit\'e dans le groupe de Bohr d'autre part, 
  qui ont introduit
 le proc\'ed\'e appel\'e aujourd'hui s\'election al\'eatoire.
  Nous allons rappeler le sens de ces termes. 
 
 Nous partirons d'une suite positive $\bw=\{w_{n}\}$ 
 et nous lui associerons une suite al\'eatoire 
 $\Lambda=\Lambda_{\bw}(\omega)$ ($\omega\in\Omega$, 
 espace de probabilit\'e) de l'une des mani\`eres \'equivalentes
que voici :

1) $\Lambda=\set {n\in \bbn}{\xi_{n}>0}$, o\`u $\xi_{n}$ 
sont des variables al\'eatoires (v.a.) ind\'ependantes, suivant des lois de Poisson de param\`etres $w_{n}\;(=\E[\xi_{n}])$.

2)  $\Lambda=\set {n\in \bbn}{\beta_{n}>0}$, o\`u $\beta_{n}$ 
sont  des variables de Bernoulli ind\'ependantes, d'esp\'erance 
$\varpi_{n}=1-e^{-w_{n}}$ $=\pr\bigbra{\beta_{n}=1}$. 

Dans les \'enonc\'es qui suivent, on peut \'ecrire $\varpi_{n}$
aussi bien que $w_{n}$. La d\'efinition 1 est utile 
dans certaines d\'emonstrations.

\begin{theoreme1}
Si $w_{n}=\osta{\frac{1}{n}}$ $(n\to\infty)$, il est presque s\^ur
que $\Lambda$ est un ensemble de Sidon, non-dense dans le groupe de Bohr $\bbb$.
\end{theoreme1}

\begin{theoreme2}
Si $\lim_{n\to\infty}nw_{n}=\infty$, il est presque s\^ur
que $\Lambda$ est un ensemble d'analyticit\'e 
(donc n'est pas  un ensemble de Sidon),
dense dans $\bbb$. 
\end{theoreme2}

Ces r\'esultats sont \`a comparer aux deux probl\`emes toujours 
ouverts :

1) (probl\`eme de ``dichotomie'') une partie de $\bbz$ est-elle
n\'ecessairement soit un ensemble de Sidon, soit
 un ensemble d'analyticit\'e ?
 
 2) (probl\`eme de l'adh\'erence dans $\bbb$ d'un ensemble de Sidon)
un ensemble de Sidon dans $\bbz$ est-il toujours non-dense
dans $\bbb$ ? On sait, \cite{tr-cluster}, que s'il existe un ensemble de Sidon dont l'adh\'erence
dans $\bbb$ contient un ouvert non vide, alors il existe
un ensemble de Sidon dense dans $\bbb$. 

$G$ \'etant un groupe ab\'elien localement compact, on d\'esigne
par $A(G)$ l'alg\`ebre des fonctions continues sur $G$ qui sont transform\'ees de Fourier de fonctions int\'egrables sur le groupe dual,
$\Gamma$. Si $A$ est un ferm\'e dans $G$, $A(E)$ d\'esigne
l'alg\`ebre des restrictions \`a $E$ des $f\in A(G)$.
Si $E$ est un compact et que $A(E)=C(E)$, espace des fonctions continues sur $E$, on dit que $E$ est un ensemble de Helson. Quand
$G$ est discret et que $A(E)=c_{0}(E)$, espace des fonctions d\'efinies
sur $E$ et tendant vers $0 $ \`a l'infini, on dit que $E$ est un ensemble de Sidon. Dans la suite, on s'int\'eresse aux cas $G=\bbt$, $\bbr$, 
$\bbz$, ou $\bbb$.

Les fonctions analytiques $F(z)$ op\`erent dans $A(G)$, c'est-\`a-dire
que $F\circ f\in A(G)$ lorsque $f\in A(G)$ et que $f$ prend ses valeurs dans le domaine de $F$. Inversement, les seules fonctions d\'efinies sur un voisinage r\'eel de $0$ qui op\`erent dans $A(G)$ sont les fonctions analytiques (nulles en $0$ si $G$ n'est pas compact).
On dit que $E$ est un ensemble d'analyticit\'e
si cela a lieu en rempla\c cant $A(G)$ par $A(E)$. Ni les ensembles
de Helson ni les ensembles de Sidon ne sont ensembles d'analyticit\'e.

$\bbb$, le compactifi\'e de Bohr de $\bbz$, est le groupe dual de
$\bbt_{d}$ ($\bbt$ discret).   
Sa topologie est la moins fine qui rende continus les caract\`eres
$\beta\rightarrow \inner{\beta}{t}$ ; une sous-base en est la famille des
ensembles 
$U(\tau,\eta,\zeta)=\set{\beta}{\abs{\inner{\beta}{\tau}-\zeta}<\eta}$,
$\tau\in \bbt$, $\eta>0$, et $\abs{\zeta}=1$.

\s
Sch\'ema des d\'emonstrations.

Si   $\bw=\bw'+\bw''$ ($w_{n}=w'_{n}+w''_{n}$), la r\'eunion de 
$\Lambda_{\bw'}$ et $\Lambda_{\bw''}$ ind\'ependants a la m\^eme
loi que $\Lambda_{\bw}$. On renforce le th\'eor\`eme 1 en augmentant
les $w_{n}$, el le th\'eor\`eme 2 en les diminuant.

\mn
{\bf Preuve du th\'eor\`eme 1. }\qua
On prend $w_{n}=\frac{\alpha}{n}$. Si $\alpha$ est assez petit
($\alpha\log 3<1$ suffit), $\Lambda$ est
un ensemble quasi-ind\'ependant, c'est-\`a-dire sans relation
lin\'eaire \`a coefficients $-1,0,1$ non tous nuls entre ses 
\'el\'ements. En g\'en\'eral, $\Lambda$ est p.s. une r\'eunion finie
d'ensembles quasi-ind\'ependants, qu'on sait \^etre un ensemble de Sidon.

La non-densit\'e r\'esulte de la proposition suivante : si $\alpha$
est assez petit ($\alpha<\alpha(\ve)$), il existe p.s. un ensemble
dense  de $t$ dans $\bbt$ tel que $[-\ve,\ve]$ contienne tous les
points de $\Lambda t$ sauf un ensemble fini.

Pour d\'emontrer la proposition, on part de la fonction triangle
d'int\'egrale $1$  et de support $[-\ve,\ve]$, soit $f$, et d'un intervalle
$I\subset \bbt$, et on \'etudie la martingale positive
\begin{equation*}
Y_{N}=\int_{I}\prod_{n=1}^{N}\bigbra{f(nt)}^{\xi_{n}}\exp\bigbra{-\frac{\alpha}{n}(f(nt)-1)}dt.
\end{equation*}
Sous la condition $\alpha<\ve^{2}$, elle converge dans 
$L^{2}(\Omega)$, et cela donne le r\'esultat voulu.
Voici les \'etapes de calcul : on pose 
$L_{N}(t)=\sum_{n=1}^{N}\frac{\cos 2\pi nt}{n}$ et on v\'erifie que
\begin{equation*}\begin{split}
\E[Y_{N}^{2}]=&\iint_{I\times I}\exp\sum_{jk\ne 0}\alpha\hat f_{j}\hat f_{k} 
L_{N}(js+kt)ds\,dt\\
\le &\prod_{jk\ne 0}\Bigbra{\iint_{I\times I}\exp \alpha\hat f_{j}\hat f_{k} 
L_{N}(js+kt)ds\,dt}^{1/p_{jk}}
\end{split}\end{equation*} avec 
$p_{jk}\inv=%
\hat f_{j}\hat f_{k}\bigbra{\sum_{jk\ne 0}\hat f_{j}\hat f_{k}}\inv$. 
En majorant $\iint_{I\times I}$ par $\iint_{\bbt\times \bbt}$,
 et $L_{N}(t)$ par $\log\frac{1}{\abs{\sin \pi t}}+C$
on obtient
\begin{equation}
\E[Y_{N}^{2}]\le \int_{\bbt}\exp \alpha\sum_{jk\ne 0}\hat f_{j}\hat f_{k}
\bigbra{\log\frac{1}{\abs{\sin \pi t}}+C}      dt,
\end{equation} 
et cette int\'egrale est born\'ee quand 
$\alpha\sum_{jk\ne 0}\hat f_{j}\hat f_{k}<1$.

\mn
{\bf Preuve du th\'eor\`eme 2.}\qua
On suppose $\lim_{n\to\infty}nw_{n}=\infty$. Posons 
$\Lambda_{N}=\Lambda\cap [1,\dots,N]$. Si $\Lambda$
\'etait un ensemble de Sidon, on aurait 
$\abs{\Lambda_{N}}=\osta{\log N}$ ($N\to \infty$), ce qui,
p.s., n'est pas le cas.

\s
{\bf Preuve que $\Lambda$ est un ensemble d'analyticit\'e.}\qua 
Elle s'inspire
de \cite{ykpm}. Quitte \`a diminuer les $w_{n}$, on suppose qu'\`a
partir d'un certain rang ils sont constants sur chaque intervalle joignant
deux multiples successifs d'une puissance de $2$ donn\'ee.
On sait que, pour tout $0<c<1$ et tout $r>0$ assez grand, il existe
$\vp\in A(\bbt)$ r\'eelle telle que $\norm{\vp}_{A(\bbt)}<r$
et $\norm{\mu e^{-i\vp}}_{PM(\bbt)}<e^{-cr}$, o\`u $\mu$ est la mesure de Haar
et $PM(\bbt)$ l'espace des pseudomesures, dual de $A(\bbt)$. 
Il s'agit, pour des $r$ arbitrairement grands, de montrer qu'il existe,
avec une probabilit\'e arbitrairement voisine de $1$, une partie finie
$\Lambda_{r}$,  une mesure $\tau_{r}$ port\'ee par $\Lambda_{r}$
et une fonction r\'eelle $\psi_{r}\in A(\bbr)$, telles que 
$\norm{\psi_{r}}_{A(\bbr)}<Cr$ et 
$\norm{\tau_{r}e^{i\psi_{r}}}_{PM(\bbr)}<C\norm{\tau_{r}}_{\mr}e^{-cr}$,
$C>0$ et $c\in (0,1)$ \'etant des constantes absolues. 

L'id\'ee est de choisir $\tau_{r}=\sum_{F}\xi_{n}\delta_{n}$, la somme \'etant prise sur un ensemble $F$ de valeurs de $n$ d\'ependant convenablement de $r$ et comportant des plages de constance pour $w_{n}$. L'\'etude d\'eterministe consiste \`a obtenir ces in\'egalit\'es
quand on remplace $\tau_{r}$ par 
$\sigma_{r}=\E[\tau_{r}]=\sum_{F}w_{n}\delta_{n}$, et l'\'etude probabiliste vise \`a montrer que 
$\norm{(\tau_{r}-\sigma_{r})e^{i\psi_{r}}}_{PM(\bbr)}$ 
est n\'egligeable au regard de $\norm{\tau_{r}}_{\mr}e^{-cr}$.

{\bf Preuve de la densit\'e dans}  $\bbt$. Elle est inspir\'ee de 
\cite{yk1},  \cite{yk2}.
Voici la proposition cl\'e : si $w_{n}=\frac{\alpha}{n}$ 
et $\alpha\abs{I}>1$ ($I$ un intervalle sur $\bbt$), il est  presque
s\^ur  que, pour tout irrationnel $t\in\bbt$, 
$\Lambda t\cap I\ne\varnothing$.    
Cette proposition admet une version
multidimensionelle, dans laquelle l'irrationnel $t\in\bbt$ est
remplac\'e par un g\'enerateur $\tb\in\bbt^{s}$, et qui entraine 
la densit\'e dans $\bbb$.

Pour \'etablir $\Lambda t\cap I\ne \varnothing$ pour tout $t$, l'id\'ee
est de se ramener \`a un ensemble fini de valeurs de $t$, quitte \`a r\'eduire $I$. Les \'etapes de la preuve sont les suivantes :

\n 1) on approche $I$
par un intervalle int\'erieur $J$, \`a distance $d$ des extr\'emit\'es de $I$

\n 2) on approche $\bbone_{J}$ dans $\lt$ par une fonction continue
$f$, $0\le f\le \bbone_{J}$

\n 3) on approche uniform\'ement $f$ par un polynome trigonom\'etrique,
de degr\'e $k$ 

\n 4) on approche $\bbt$ par un ouvert 
$G=\set{t}{\forall j\in\{1,2,\dots,k\}\;\abs{\sin \pi jt}>\delta}$ et on
v\'erifie que 
pour un $\beta$ arbitrairement proche de $\alpha\abs{I}$
et pour un $C=C(\beta,\delta )$ fini on a \\ \centerline{
$\int_{G }  \pr\bigbra{\Lambda_{N}t\cap J=\varnothing}dt<CN^{-\beta}$}

\n 5) on choisit $M=\frac{N}{d}$ et on v\'erifie qu'il existe un $\vt\in\bbt$
tel que
\begin{equation*}%\label{}
 \sum_{\vt+mM\inv\in G}\pr\bigbra{%
 \Lambda_{N}(\vt+mM\inv)\cap J=\varnothing
 } < CMN^{-\beta}=o(1) (N\to\infty)
\end{equation*}
6) on approche tout $t\in G$ donn\'e par un $\vt+mM\inv\in G$ 
qui en est
\`a distance $<d$, d'o\`u $\pr (\Lambda t\cap I=\varnothing)=0$.

 \sn
 {\bf Remarques sur les suites d'entiers dans le groupe de Bohr} $\bbb$.
 
 Diff\'erents crit\`eres de densit\'e et des exemples sont donn\'es dans 
 \cite{yk2}.
 
 Les ensembles $I_{0}$ de Hartman et Ryll-Nardzewski
forment une classe particuli\`ere d'ensembles de Sidon, et leur
adh\'erence dans $\bbb$ est un ensemble de Helson. Il a \'et\'e
 conjectur\'e que tout ensemble de Sidon soit une r\'eunion finie d'ensembles $I_{0}$. 
 Il s'ensuivrait que l'adh\'erence d'un ensemble de Sidon
soit un ensemble de Helson dans $\bbb$, \cite{tr-comp}.

En compl\'ement du probl\`eme 2, on peut donc se demander
si l'adh\'erence dans $\bbb$ d'un ensemble de Sidon
est n\'ecessairement un ensemble de Helson.
Est-ce n\'ecessairement un ensemble de mesure de Haar nulle ?

En ce qui concerne les entiers al\'eatoires, le th\'eor\`eme 1 r\'epond \`a
la question de non-densit\'e, mais laisse ouverte la question de la mesure de Haar de l'adh\'erence de la suite $\Lambda$.

\renewcommand{\refname}{\normalfont\itshape{R\'ef\'erences}}

\mn
Jean-Pierre Kahane, $<$jean-pierre.kahane@math.u-psud.fr$>$, \\
Yitzhak Katznelson, $<$katznel@math.stanford.edu$>$.
\end{document}